\documentclass[12pt]{amsart}
\usepackage{amscd,amsmath,amssymb,amsfonts}
\usepackage[cmtip, all]{xy}
\theoremstyle{plain}
\newtheorem{thm}{Theorem}
\newtheorem{lem}[thm]{Lemma}

\newtheorem{prop}[thm]{Proposition}

\theoremstyle{definition}

\newtheorem{rmk}[thm]{Remark}

\numberwithin{thm}{section}
\numberwithin{equation}{section}

\newcommand{\ml}[2]{\begin{multline}\label{#1}#2 \end{multline}}
\newcommand{\ga}[2]{\begin{gather}\label{#1}#2 \end{gather}}

\newcommand{\sF}{{\mathcal F}}

\newcommand{\sI}{{\mathcal I}}
\newcommand{\sJ}{{\mathcal J}}
\newcommand{\sK}{{\mathcal K}}
\newcommand{\sL}{{\mathcal L}}

\newcommand{\sO}{{\mathcal O}}

\newcommand{\sV}{{\mathcal V}}

\newcommand{\A}{{\mathbb A}}

\newcommand{\C}{{\mathbb C}}

\renewcommand{\H}{{\mathbb H}}

\newcommand{\N}{{\mathbb N}}
\renewcommand{\P}{{\mathbb P}}

\begin{document}

\title[Artin's vanishing]{Variation on Artin's vanishing theorem
}
\author{H\'el\`ene Esnault}
\address{
Universit\"at Duisburg-Essen, FB6, Mathematik, 45117 Essen, Germany}
\email{esnault@uni-essen.de}
\thanks{Partially supported by the 
DFG-Schwerpunkt ``Komplexe Mannigfaltigkeiten'' and by the DFG Leibniz Preis.}

\date{June 20, 2004}
\dedicatory{\`A Michael Artin, \`a l'occasion de  son 70-i\`eme anniversaire}
\begin{abstract}
We give a proof of Artin's vanishing theorem in characteristic zero, based on 
Deligne's Riemann-Hilbert correspondence. 
\end{abstract}
\maketitle
\begin{quote}

\end{quote}

\section{Introduction}
In \cite{Ar}, Corollaire 3.5, M. Artin shows that if $X$ is  a affine 
variety over a separably closed field $k$, and if
$\sF$ is a constructible sheaf on it, then \'etale cohomology 
$H^m_{\acute{e}t}(X, \sF)$ is vanishing  for $m\ge {\rm dim}(X)+1$.   He reduces the proof to dimension 1 by applying base change for proper morphisms. The purpose of this note is to give an analytic proof of Artin's 
famous vanishing theorem in the analytic category, that is when $k$ is the field of complex numbers, $\sF$ is a constructible sheaf of complex vector spaces
with possibly infinite monodromies on strata, and the cohomology is the analytic one. Surely this is not necessary to have an analytic proof to apply Artin's vanishing theorem analytically, but this is perhaps 
of interest to know it exists.   
Indeed, we give two proofs in the framework of Deligne's Riemann-Hilbert correspondence \cite{DePSR}. The first one is based on a splitting property proven in \cite{E}, Proposition 1.2. The second one is a slightly different
 form of this
and is due to  Pierre Deligne. He did not  write it down in his Lectures Notes 
\cite{DePSR}, 
but communicated to us in a discussion on the first proof. It is based on 
\cite{DeRD}, Proposition 5, p. 412. In both cases, it allows to compute the analytic cohomology of the extension by $0$ of a local system on a Zariski open set as the hypercohomology of a suitable algebraic 
de Rham complex.
\\ \ \\
{\it Acknowledgement:} We thank Pierre Deligne for explaining to us his proof.
 The analytic proof of the vanishing theorem was 
explained in  a course held in the spring  2003 at the University of Essen. We 
thank the students for their active participation.

\section{ Artin's vanishing in complex geometry}
\begin{thm} \label{thm2.1} Let $X$
be an affine variety defined over the field of complex numbers, and
$\sF$ be a constructible sheaf. Then  $H^m(X_{{\rm an}}, \sF)=0 $ 
for $m> {\rm dim}(X)$. 
\end{thm} 
\begin{proof} 
If $Y\subset X$ is the support of $\sF$, then $H^m(X_{{\rm an}}, \sF)=H^m(Y_{{\rm an}}, \sF)$ and 
$Y$ is affine. Thus we may assume that the support of $\sF$ is $X$. Let $\pi: X\to \A^n$ be a Noether normalization, with $n={\rm dim} (X)$. Then 
$H^m(X_{{\rm an}},\sF)=H^m(\A^n_{{\rm an}}, \pi_*\sF)$, and $\pi_*\sF$ is a constructible sheaf with support $\A^n$. Thus we may assume $X=\A^n$. Let $j: \emptyset \neq U
\hookrightarrow X$ be a Zariski open set so that $\sF|_{U_{{\rm an}}}=\sV$ is a local
 system. Thus $j_!\sV\subset \sF$ and the quotient $\sF/j_!\sV$ is
a constructible sheaf with support in 
dimension $< n$. Thus inducting on ${\rm dim}(X)$, we may assume 
$\sF=j_!\sV$, as for ${\rm dim }(X)=0$ the theorem is trivial. 
Let $Y:=X\setminus U$. If $Y$ has codimension $\ge 2$, then $\sV$ is the constant sheaf $\C^r$, where $r$ is the rank of $\sV$. From the exact sequence
$0\to j_!\C^r \to \C^r \to \C^r|_Y\to 0$, and the induction, we reduce 
the theorem to the vanishing of $H^m(\A^n_{{\rm an}}, \C)=0$ for $m> n$, which is trivial. 
Else $Y$ is a divisor. Then there is an algebraic bundle $V$ on $X$ (which has  to be trivial since $X =\A^n$)  with an algebraic connection $\nabla: 
 V \to \Omega^1_X(mY)\otimes V$ with meromorphic poles along $Y$, such that 
$(V_{{\rm an}}|_U)^\nabla=\sV$. 
We choose for $(V,\nabla)$ any extension to $X$ of the 
unique regular singular connection on $X\setminus Y$ (\cite{DePSR}, Th\'eor\`eme 5.9). In particular, it is regular singular ``at the $\infty$ of $X$''. 
 Let $\sI$ be the ideal sheaf of $Y$. For $N \in \N, N>n$, one defines the complex
\ml{2.1}{\sK_N: \sI^N \otimes V\to 
\sI^{N-1}\otimes V\otimes \Omega^1_X(mY) \to \ldots\\ \ldots \to 
\sI^{N-n} \otimes V \otimes \Omega^n_X(nmY) } 
\begin{lem} \label{lem2.2}
For  $N\in \N$ large,  $H^m(X_{{\rm an}}, j_!\sV)$ is a direct summand of
$\H^m(X, \sK_N)$. 
\end{lem}
\begin{proof} (See \cite{E}, Proposition 1.2). 
Let $\pi:\P \to \P^n$ be a birational projective morphism, with $\P$ smooth, 
so that $\bar{Z}+H$ is a normal crossing divisor. Here we denote by
 $Z$ the inverse image $\pi^{-1}(Y)$, by $\bar{Z}$ its Zariski closure
in $\P$ and by  $H$ the inverse image $\pi^{-1}(\P^n\setminus \A^n)$. 
We also set $\lambda: \A^n\to \P^n, \tilde{\lambda}: 
\pi^{-1}(\A^n)=\P\setminus H\to \P, j: U\to \A^n, \tilde{j}: U\to 
\P\setminus H$. 
Let $(\bar{V},\bar{\nabla})$ be an extension to $\P$ 
of $(V,\nabla)$ on $U$, so that $\bar{V}$ is locally free and
$\bar{\nabla}$ has logarithmic poles along $\bar{Z}+H$. Then for $M$ large enough, $R\tilde{\lambda}_*\tilde{j}_!V\to 
(\Omega^\bullet_{\P}(\log (\bar{Z}+H))(-M\bar{Z})\otimes \bar{V})_{{\rm an}}
\to (\Omega^\bullet_{\P}(\log \bar{Z})(*H)(-M\bar{Z})\otimes \bar{V})_{{\rm an}}$ are quasi-isomorphisms
(\cite{DePSR}, Corollaire 6.10, and via duality \cite{EV}, (2.9), (2.11)), thus by GAGA \cite{Se}, it induces an isomorphism
$H^m(\P\setminus H, \tilde{j}_!\sV)(=H^m(\A^n, j_!\sV))\xrightarrow{\cong}
\H^m(\P \setminus H, \Omega^\bullet_{\P\setminus H}(\log Z)(-M Z)
\otimes \bar{V})$. Taking now any  coherent extension $V'$ to $\P^n$ of $V$, with $\sI'$ being the ideal sheaf of the Zariski closure $Y'$  of $Y$ in $\P^n$, 
one has $j_*\sK_N=(\sI')^{N-\bullet}\otimes V' \otimes 
\Omega^\bullet_{\P^n}(\bullet m  Y')(*(\P^n \setminus \A^n))$. By GAGA again, one has
$H^m(\A^n, \sK_N)=
\H^m(\P^n, (\sI')^{N-\bullet}\otimes V' \otimes 
\Omega^\bullet_{\P^n}(\bullet m  Y')(*(\P^n \setminus \A^n)))=
\H^m(\P^n, ((\sI')^{N-\bullet}\otimes V' \otimes 
\Omega^\bullet_{\P^n}(\bullet m  Y')(*(\P^n \setminus \A^n)))_{{\rm an}}).$
The latter group receives $H^m(\A^n_{{\rm an}}, j_!V)$ while functoriality implies  the existence of $$\pi^{-1}: \sK_N\to 
R\pi_*(\Omega^\bullet_{\P\setminus H}(\log Z)(-M Z)\otimes \bar{V}),$$ which is an isomorphism on $U$. This finishes the proof.
\end{proof}
 As $X$ is affine, one has
$\H^m(X, \sK_N)= H^m(\Gamma(X,\sK_N))=0$ for $m>{\rm dim}(X)$. Using 
lemma \ref{lem2.2} this concludes the proof. 
\end{proof} 
\begin{rmk}
Of course in Lemma \ref{lem2.2}, we could replace $X=\A^n$ by any smooth 
variety $X$.
\end{rmk}
We now reproduce a communication by P. Deligne, which completes Lemma \ref{lem2.2}. Instead of considering \eqref{2.1} for a fixed $N$, which leads to a splitting of $j_!$, we consider the inverse system of such. More generally, let $\sV$ be a local system on a smooth Zariski open $\emptyset \neq U$
of a complex variety $X$, and let  $\sK$ be
a bounded complex of coherent sheaves,  extending 
$V_0\otimes \Omega^\bullet_{U}$ to $X$, 
where $(V_0,\nabla_0)$ is the unique algebraic bundle with a regular singular connection on $U$ with underlying $\sV$, \cite{DePSR}, loc. cit. 
Let $\sI$ be a sheaf of ideals with supports $Y:=X\setminus U$, and let 
$(\sL_N)_N$ be the projective system of complexes $\sL_N$ defined by 
\ga{2.2}{(\sL_N)^p= \sI^{N-p}\sK^p:={\rm Im} (\sI^{N-p}\otimes_{\sO_X} \sK^p \to \sK^p).}
For example, for $\sK=\Omega_X^\bullet(\bullet  mY)\otimes V$ as in the proof of Theorem \ref{thm2.1}, then $\sL_N$ is nearly equal to $\sK_N$  of \eqref{2.1} (nearly as we haven't assumed $V$ to be locally free, thus the tensor product is not necessarily equal to the product in \eqref{2.2}.) We denote by $j:U\to X$ the open embedding. 
\begin{prop}[P. Deligne] One has
$$H^m(X_{{\rm an}}, j_!\sV)=\varprojlim_N \H^m(X, \sL_N).$$
\end{prop}
\begin{proof}(P. Deligne)  Let 
$\sigma: \tilde{X}\to X$ be a projective birational morphism with $\sigma|_U={\rm id}$, with $\tilde{X}$ smooth and $\sigma^{-1}(Y)=:Z$ a normal crossing divisor. Let us set $\sJ=\sigma^*\sI$ and choose $(V, \nabla)$ 
an extension to $\tilde{X}$ of $(V_0,\nabla_0)$ on $U$. By \cite{DeRD}, Proposition 5,  
the projective system $R^a\sigma_*
(\sJ^N \otimes V\otimes \Omega^p_{\tilde{X}}(\log Z))$ is essentially constant of value $0$ for $a>0$. This means that for $N\ge 0$ given, there is a $N'>N$ so that $
R^a\sigma_*
(\sJ^{N'} \otimes V\otimes \Omega^p_{\tilde{X}}(\log Z))
\xrightarrow{0} R^a\sigma_*
(\sJ^N \otimes V\otimes \Omega^p_{\tilde{X}}(\log Z))$. 
  Moreover, for  all $N$, there are $N_i, i=1,2$ with 
$(\sL_{N_2})^p\to \sigma_* (
V\otimes \sJ^{N_1}\otimes \Omega^p_{\tilde{X}}(\log Z) )\to (\sL_N)^p$. 
This shows that in the proposition, we may assume that $X$
is smooth, $Z=X\setminus U$ is a normal crossing divisor, $\sK^\bullet= 
V\otimes \Omega^\bullet_X(\log Z)$, where $(V,\nabla)$ extends $(V_0,\nabla_0)$ with $V$ locally free and $\nabla$ with logarithmic poles. 
Again for all $N$, there there are $N_i, i=1,2$ with 
$\sL_{N_2}\to 
V\otimes \sO_X(-N_1Z)\otimes \Omega^\bullet_{X}(\log Z) )\to \sL_N$.
Thus $\varprojlim_N \H^m(X, \sL_N)=\varprojlim_N 
\H^m(X, V\otimes \sO_X(-N Z)\otimes \Omega^\bullet_{X}(\log Z) )$. 
Taking $\bar{X}\supset X$ a good compactification with 
$\infty=\bar{X}\setminus X$ a normal crossing divisor, and 
with $(\bar{V}, \bar{\nabla})$ an extension to $\bar{X}$  of $(V,\nabla)$ with $\bar{V}$ locally free and $\bar{\nabla}$ with logarithmic poles,  
one has   
$\H^m(X, V\otimes \sO_X(-N Z)\otimes \Omega^\bullet_{X}(\log Z) )=
\H^m(\bar{X}, \bar{V}\otimes \sO_{\tilde{X}}(-N \bar{Z})
\otimes \Omega^\bullet_{\bar{X}}(\log \bar{Z})(*\infty))$. Here
 $\bar{Z}\subset \bar{X}$ is the Zariski closure of $Z\subset X$. 
 Again by GAGA \cite{Se}, 
\ml{}{
\H^m(\bar{X}, \bar{V}\otimes \sO_{\bar{X}}(-N \bar{Z})
\otimes \Omega^\bullet_{\tilde{X}}(\log \bar{Z})(*\infty))= \notag \\
\H^m(X_{{\rm an}}, (V\otimes \sO_X(-N Z)\otimes 
\Omega^\bullet_{X}(\log Z))_{{\rm an}} )} and the latter for $N$ large enough is equal to $H^m(X_{{\rm an}}, j_!\sV)$ by 
\cite{DePSR}, Corollaire 6.10, and via duality \cite{EV}, (2.9), (2.11).

\end{proof}
\bibliographystyle{plain}

\renewcommand\refname{References}

\end{document}